\documentclass[12pt]{amsart}

\usepackage{amsmath, amsthm, amssymb, amsfonts, verbatim}
\usepackage{epsfig}

\def\0{{\bf 0}}

\def\R{{\mathbb R}}
\def\Z{{\mathbb Z}}

\def\N{{\mathbb N}}

\newtheorem{theorem}{Theorem}[section]

\theoremstyle{definition}
\newtheorem{definition}[theorem]{Definition}

\theoremstyle{remark}

\numberwithin{equation}{section}

\begin{document}

\newcommand{\Obsolete}[1]{
    }
%
\newcommand{\nwc}{\newcommand}
\nwc{\qref}[1]{(\ref{#1})}
\nwc{\bsa}{\beta_{\sigma,\alpha}}
\nwc{\bpk}{\hat\beta_{+,k}}
\nwc{\bmk}{\hat\beta_{-,k}}
\nwc{\bk}[1]{\hat\beta_{#1}}
\nwc{\KK}{{\mathcal K}}
\nwc{\Ks}{{\mathcal K}_{\sigma}}
\nwc{\ps}{{p_{\rm S}}}
\nwc{\lap}{\Delta}
\nwc{\grad}{\nabla}
\nwc{\bkp}{\beta_{+,k}}
\nwc{\bkm}{\beta_{-,k}}
\renewcommand{\div}{\nabla\cdot}
\nwc{\nddu}{n\cdot(\Delta-\nabla\nabla\cdot)u}
\title[The Laplace-Leray commutator in domains with corners]{On optimal estimates for the Laplace-Leray commutator in planar
domains with corners}

\author{Elaine Cozzi}
\address{Department of Mathematical Sciences, Carnegie Mellon University}
\email{ecozzi@andrew.cmu.edu}
\author{Robert L. Pego}
\address{Department of Mathematical Sciences, Carnegie Mellon University}
\email{rpego@andrew.cmu.edu}
\thanks{This material is based upon work supported by the National Science Foundation under Grant Nos.\  DMS06-04420 and DMS09-05723 and partially supported by the Center for Nonlinear Analysis (CNA) under National Science Foundation Grant No.\ DMS06-35983.}
\subjclass{Primary 35} 

\begin{abstract}
For smooth domains, 
Liu et al.~(Comm. Pure Appl. Math. 60: 1443-1487, 2007) used optimal estimates for the commutator of the Laplacian and 
the Leray projection operator 
to establish well-posedness of an extended Navier-Stokes dynamics.
In their work, the pressure is not determined by incompressibility, but rather by a certain formula involving the Laplace-Leray commutator.
A key estimate of Liu et al.\ controls the commutator
strictly by the Laplacian in $L^2$ norm at leading order.
In this paper we show that this strict control fails in a 
large family of bounded planar domains with corners.
However, when the domain is an infinite cone,
we find that strict control may be recovered in 
certain power-law weighted norms.
\end{abstract}
\maketitle
\section{Introduction}\label{Introduction}
In this paper, we study estimates for $[\Delta,P]=\Delta P-P\Delta$, the commutator of the Laplacian and the Leray projection operator, in planar domains with corners.  
In a bounded domain $\Omega\subset\R^N$, the Leray
projection operator $P$ is defined as follows:  Given any $a\in
L^2(\Omega,\R^N)$, there exists a unique $q\in H^1(\Omega)$ with
$\int_{\Omega} q=0$ and such that $Pa := a + \nabla q$ satisfies  
\begin{equation}\label{helmholtz}
0=\langle Pa,\nabla \phi \rangle = \langle a+\nabla q, \nabla \phi \rangle 
\end{equation}
for all $\phi\in H^1(\Omega)$.  In \cite{LLP}, 
Liu et al.\ proved the following $L^2$-estimate for the commutator of the Leray
projection operator and the Laplacian.
\begin{theorem}\label{mainLLP}
Let $\Omega$ be a connected, bounded domain in $\R^N$, $N\geq 2$, with
$C^3$ boundary.  For any $\beta>\frac12$, there exists $C\geq 0$ such
that for all vector fields $u\in H^2\cap H^1_0(\Omega,\R^N)$,
\begin{equation}\label{stabilityest}
\int_{\Omega} |[\Delta,P] u|^2 
\leq 
\beta
\int_{\Omega} |\Delta u |^2 + C \int_{\Omega} |\nabla u|^2. 
\end{equation} 
\end{theorem}   

Theorem \ref{mainLLP} has significant applications to the
Navier-Stokes equations.  We recall that on a bounded domain $\Omega$
in $\R^N$ for $N\geq 2$, the Navier-Stokes equations modeling
incompressible viscous fluid flow with no-slip boundary conditions are
given by 
\begin{align*}
 \begin{matrix}
   (NS) & \left\{
      \begin{matrix}
         \partial_t u + u \cdot \nabla u  + \nabla p = \nu \Delta u + f \\
             \nabla\cdot u = 0 \\
             u|_{\Gamma}=0,
         \end{matrix}
        \right.
    \end{matrix}
\end{align*}
where $\Gamma=\partial\Omega$, $u$ denotes the
velocity of the fluid, $p$ denotes the pressure,
and $\nu$ represents the viscosity.  In \cite{LLP}, the authors
consider strong solutions to ($NS$) with constant $\nu>0$, 
and show that the pressure satisfies 
\begin{equation}\label{E.p}
\nabla p = (I-P)(f-u\cdot \nabla u)+ \nu [\Delta,P] u.
\end{equation}
For such solutions they prove the unconditional stability and
convergence of a simple time discretization scheme which decouples the
updates of velocity and pressure.  
The decoupling of these variables is significant in that it eliminates
the need for an inf-sup condition which is often necessary to prove
the stability in finite-element schemes.  
A critical ingredient in the proof of stability in \cite{LLP} 
is that by invoking Theorem \ref{mainLLP} with $\beta<1$, 
one can strictly control the pressure gradient by the viscosity term
plus lower-order terms.  
As a result, Liu et al.\ establish the well-posedness of an {extended}
Navier-Stokes dynamics in which the pressure $p$ is always determined 
by the formula (\ref{E.p}) and the zero-divergence condition
is dropped in general.
We refer the reader to \cite{LLP} for further details and discussion.

\Obsolete{In an attempt to motivate the significance of such an
estimate, in this section we briefly describe the connection between
the commutator and the unconstrained formulation of ($NS$) introduced
in \cite{LLP}.  We refer the reader to \cite{LLP} for a more thorough
discussion of this connection.  We also refer the reader to \cite{LLP}
for the details of the time discretization scheme used to approximate
solutions to this formulation.

To establish their unconstrained formulation of ($NS$), Liu, Liu and
Pego first use the property $\nabla\cdot u = 0$ to  rewrite ($NS$) in
the form 
\begin{equation}\label{NS1}
\partial_t u + P(u\cdot\nabla u -f-\nu\Delta u) = \nu\nabla \nabla\cdot u.
\end{equation}
The authors then use the equality $\nabla\nabla\cdot u=\Delta(I-P)u$
(Lemma 1.1 of \cite{LLP}), to see that
\begin{equation}\label{commform1}
[\Delta,P] u = (I-P)\Delta u -\nabla \nabla\cdot u 
= (I-P) ( \Delta-\nabla\nabla\cdot) u,
\end{equation}
and they rewrite ($\ref{NS1}$) as 
\begin{equation}\label{NS2}
\partial_t u + P(u\cdot\nabla u -f) + \nu[\Delta, P]u= \nu\Delta u.
\end{equation}
Comparing ($NS$) with ($\ref{NS2}$), we observe that the pressure
gradient $\nabla p$ in ($NS$) takes the form
\begin{equation}\label{pressure}
\nabla p = (I-P)(f-u\cdot \nabla u)+ \nu [\Delta,P] u.
\end{equation}  
In \cite{LLP}, the authors show that the $L^2$-norm of the second term
of ($\ref{pressure}$) is dominated strictly by the $L^2$-norm of the
viscosity term in ($NS$).}

Theorem \ref{mainLLP} assumes that the boundary $\Gamma$ of $\Omega$
is $C^3$.  One would like to weaken this assumption to allow, for
example, sharp corners on $\Gamma$. In this paper, we show that such an improvement is not possible.  We let $\mathcal{K}_{\sigma}$ denote an infinite cone centered at the
origin, taking the form  
\begin{equation}\label{E.dom}
{\mathcal{K}}_{\sigma}=\{ (x_1,x_2)\in\R^2:0<r<\infty, 0<\theta < \sigma\},
\end{equation}
where $r$ and $\theta$ denote the polar coordinates of $(x_1,x_2)$ and
$\sigma\in(0,2\pi)$.  We consider bounded domains $\Omega\subset\R^2$
satisfying the following property: there is a neighborhood $U$ of $0$ 
such that $U\cap \tilde\Omega=U\cap\mathcal{K}_{\sigma}$
for some rotated translate $\tilde\Omega=R(\Omega-x_0)$ of $\Omega$
and for some $\sigma\ne\pi$.
In this case we call $\Omega$ a {\em bounded domain with a straight corner}.  
We claim that Theorem \ref{mainLLP} fails on any such domain.  
\begin{theorem}\label{boundedcase}
Let $\Omega$ in $\R^2$ be a bounded domain with a straight corner.  Then for every
$\beta<1$ and for every $C\in\R$, there is a vector field
$u\in H^2\cap H^1_0(\Omega,\R^2)$ satisfying 
\begin{equation}
\int_{\Omega} |[\Delta,P] u|^2 > \beta \int_{\Omega} |\Delta u |^2 + C
\int_{\Omega} |\nabla u|^2. 
\end{equation} 
\end{theorem}

One may suspect that the reason $\beta<1$ is not possible in general
has something to do with the lack of $H^2$ regularity for the Stokes
operator in domains with reentrant corners. One known way of
dealing with this situation involves using {\em weighted} Sobolev spaces.  
In a recent paper of Rostamian and Soane \cite{RS}, the authors
reformulate the time discretization scheme of \cite{LLP} in non-convex
polygonal domains using such {weighted} spaces.  
While the authors do not prove convergence of their scheme, they do
give numerical evidence suggesting that this scheme converges to the
correct solution.

We are motivated by \cite{RS} 
and elliptic regularity theory with weights \cite{KMR} 
to allow for corners on $\Gamma$ and
look to prove an optimal estimate similar to ($\ref{stabilityest}$) in
a weighted $L^2$-space.  For the most part, we study
conical domains of the form in (\ref{E.dom}).
The weighted spaces considered in \cite{KMR} are defined as follows.
\begin{definition}\label{weighted}
For an integer $l\geq 0$ and a real number
$\alpha$, we define the space $V^l_{2,\alpha}({\mathcal{K}}_{\sigma})$ to be the
closure of $C^{\infty}_c(\overline{\mathcal{K}}_{\sigma}\backslash \{0\})$ with respect
to the (scale-invariant) norm
\begin{equation}
\|f\|_{V^l_{2,\alpha}} = \left(\int_{{\mathcal{K}}_{\sigma}} { \sum_{|\rho|\leq l}
r^{2(\alpha-l+|\rho|)} |D^{\rho}_x f |^2} r\, dr\,d\theta
\right)^{\frac{1}{2}} <\infty.
\end{equation}
\end{definition}
\noindent We refer the reader to \cite{KMR} for a more thorough
discussion of weighted Sobolev spaces in an infinite cone.  

Before we state the main theorem, we must define the Leray projection operator
on unbounded domains.  This definition differs from that given in
($\ref{helmholtz}$), because if $\Omega$ is unbounded, then $\nabla
H^1(\Omega)$ is not closed in $L^2(\Omega)$.  To remedy this, we fix a bounded
domain $B\subset\Omega\subset\R^N$, and we define the space 
\begin{equation}\label{projunbounded}
Y = \left \{ q\in L^2_{loc} (\Omega): \nabla q \in L^2(\Omega,\R^N) \text{ and } \int_{B} q =0 \right \}.
\end{equation}   
Then $Y$ is a Hilbert space with norm $\|q\|^2_{Y}=\int_{\Omega} |\nabla q|^2$, and the space $\nabla Y$ is closed in $L^2(\Omega,\R^N)$.  We define the Leray projection operator $P$ as in ($\ref{helmholtz}$), except that we assume $q$ is in $Y(\Omega)$ instead of $H^1(\Omega)$.  Further discussion of the Leray projection operator on unbounded domains can be found in \cite{S}.

We remark that if $\Omega$ is Lipschitz, $C_c^{\infty}(\overline{\Omega})$ is
dense in $Y$.  The proof of this fact is similar to the proof for $\Omega=\R_+^N$ 
indicated in \cite{LLP}, based on the case $\Omega=\R^N$ treated in \cite[Lemma~2.5.4]{S}.  

We are now prepared to state the main theorem. 
\begin{theorem}\label{main}
Suppose $\sigma\in(0,2\pi)$ and let ${\mathcal{K}}_{\sigma}$ be an infinite
planar cone as in (\ref{E.dom}).
Let $\alpha\ne1$.  Then the following estimate holds for all $u\in
C^{\infty}_c(\overline{\mathcal{K}}_{\sigma}\backslash \{ 0\},\R^2)$:
\begin{equation}\label{estonweighted}
\int_{{\mathcal{K}}_{\sigma}} {r^{2\alpha} |[\Delta,P] u|^2 r} \,dr\,d\theta 
\leq \bsa\int_{{\mathcal{K}}_{\sigma}} { r^{2\alpha} |\Delta u|^2 r }\,dr\,d\theta ,
\end{equation}
where 
\begin{equation*}
\bsa= \sup_{k>0} \max 
\left\{ \bpk,\bmk \right\},
\end{equation*}
with 
\begin{equation}\label{E.bpmk}
\bk{\pm,k}= \frac{k^2+\alpha^2}{2k^2(1-\alpha)}(1-e^{-2k\sigma})  
\Re \left\{  \frac{(1-\alpha+ik)(1 \pm e^{-(k+i\alpha-2i)\sigma})}
{(1 \pm e^{-(k-i\alpha)\sigma})(1- e^{-2(k+i\alpha-i)\sigma}) } \right\}.
\end{equation}
Moreover, $\bsa$ is the smallest constant satisfying
$(\ref{estonweighted})$ for every 
$u\in C^{\infty}_c(\overline{\mathcal{K}}_{\sigma}\backslash \{ 0\},\R^2)$.
\end{theorem}

We will prove Theorem \ref{main} in Sections \ref{Preliminaries}
and \ref{solveforbeta}.  In Section \ref{bounded}, we show that
Theorem \ref{main} implies Theorem \ref{boundedcase}.

The expressions in (\ref{E.bpmk}) are sufficiently complicated that
it is difficult to characterize exactly when $\bsa<1$ holds.
We will make a few observations, however, and provide numerical evidence 
which suggests that for all $\sigma\in(0,2\pi)$ {\em except for one value}
$\sigma=\sigma_c \approx 1.4303\pi$, we have $\bsa<1$ for $\alpha$ in some interval
just to the left or right of $\alpha=0$. 

First, note that as $k\to\infty$ we have $\bk{\pm,k}\to\frac12$.
For $\alpha=0$ we compute that
\begin{equation}\label{alpha0}
\bk{\pm,k}= \frac12
\frac{\cosh^2 k\sigma - \cos^2 \sigma
\mp \cosh k\sigma\sin^2 \sigma \mp k\sin \sigma\cos \sigma \sinh k\sigma}
{\cosh^2 k\sigma-\cos^2 \sigma},
\end{equation}
from which we see 
that if $\sigma=\pi$, then $\bk{\pm,k}\equiv\frac12$, hence $\beta_{\pi,0}=\frac12$.
This half-space estimate (\ref{estonweighted}) with constant weight was already proved 
in \cite{LLP}, and explains why the condition $\beta>\frac12$ is essentially optimal 
in Theorem~\ref{mainLLP}. 
Note that due to the dilation invariance of the domain, no lower-order term
such as that in (\ref{stabilityest}) should appear in the half-space case, since it would scale differently 
under dilation.

Whenever $\pi\ne\sigma\in(0,2\pi)$ however, we have
$\bk{-,0}=1,$ $\bk{+,0}=0$.  
Thus, whenever the weight is constant ($\alpha=0$) and the cone has a corner
($\sigma\ne\pi$) we conclude that the optimal constant $\beta_{\sigma,0}\ge1$.
Our proof of Theorem~\ref{boundedcase} relies on this fact.

It is easy to approximate $\bsa$ numerically. 
For a number of values of the cone angle $\sigma$, in Figure 1
we plot $\log_{10}\bsa$ vs.~$\alpha$ for $\alpha\in[-1,1]$.
Spikes appear in many of these graphs, 
providing evidence of singularities where presumably $\bsa=+\infty$. 
After closer examination, these graphs suggest that:
\begin{itemize}
\item $\bsa=1$ whenever $\alpha=0$ and $\sigma\ne\pi$.
\item $\bsa<1$ for small $\alpha>0$ when $0<\sigma<\pi$ or $\sigma_c<\sigma<2\pi$.
\item $\bsa<1$ for small $\alpha<0$ when $\pi<\sigma<\sigma_c$. 
\end{itemize}
The number $\sigma_c \approx 1.4303\pi$ satisfying 
$\sigma_c\cot\sigma_c-1=0$ appears to be a critical value of $\sigma$
where the minimum of $\bsa$ occurs at $\alpha=0$, and the minimum value is 1.
To see this, we observe that numerical evidence indicates that for $\sigma$
near the critical value and for $\alpha$ near $0$, $\bsa$ is achieved at $k=0$.
We therefore take the limit as $k$ approaches $0$ of $\bpk$ and $\bmk$, which
yields the formulas for $\hat{\beta}_{+,0}$ and $\hat{\beta}_{-,0}$ given in
(\ref{blowup}) and (\ref{blowup2}).  Numerical evidence again shows that for
$\sigma$ in a neighborhood of the critical value and for $\alpha$ near $0$,
$\hat{\beta}_{+,0}>\hat{\beta}_{-,0}$.  Using Maple to differentiate $\hat{\beta}_{+,0}$
with respect to $\alpha$, and evaluating the derivative at $\alpha=0$, we find
that 
\begin{equation}
\partial_{\alpha} \hat{\beta}_{+,0}|_{\alpha=0}=\sigma\cot\sigma - 1.
\end{equation}


These numerical results also suggest that $\bsa<1$ for convex
cones, uniformly in $\sigma$ for positive $\alpha$ in a fixed interval.
So we may conjecture that for a bounded polygonal domain $\Omega$
that is convex, say, an estimate of the form
\begin{equation}\label{conjecture}
\int_{\Omega} {r^{2\alpha} |[\Delta,P] u|^2 } 
\leq \int_{\Omega} 
r^{2\alpha} \left( \beta   |\Delta u|^2  
+ C |\grad u|^2 \right) 
\end{equation}
will hold for some $\beta<1$ and $C$ independent of $u$ in a suitable space
of functions vanishing on $\partial\Omega$,
provided $\alpha$ is small and positive. Here $r=r(x)$ would be the 
distance from $x\in\Omega$ to the nearest corner on $\Gamma$.  The lower order term on the right hand side of (\ref{conjecture}) comes from the definition of the $V^2_{2,\alpha}$ norm on a {\em bounded} domain (see \cite{KMR}), given by
\begin{equation*}
\| u \|_{V^2_{2,\alpha}(\Omega)} = \left(\int_{\Omega} r^{2\alpha} \sum_{|\rho| \leq 2} |D_x^{\rho} u|^2 \, dx\right)^{\frac{1}{2}}.
\end{equation*}  
We do not include the term $\| r^{\alpha} u  \|^2_{L^2}$ on the right hand side of (\ref{conjecture}), because it can be controlled by first order partial derivatives using a Hardy inequality (see \cite{KMR}, Chapter 7 for details).    

However, we have no proof of (\ref{conjecture}) at this time.
 {
\begin{figure}
\label{F.beta}
\centerline{\epsfysize=5in{\epsffile{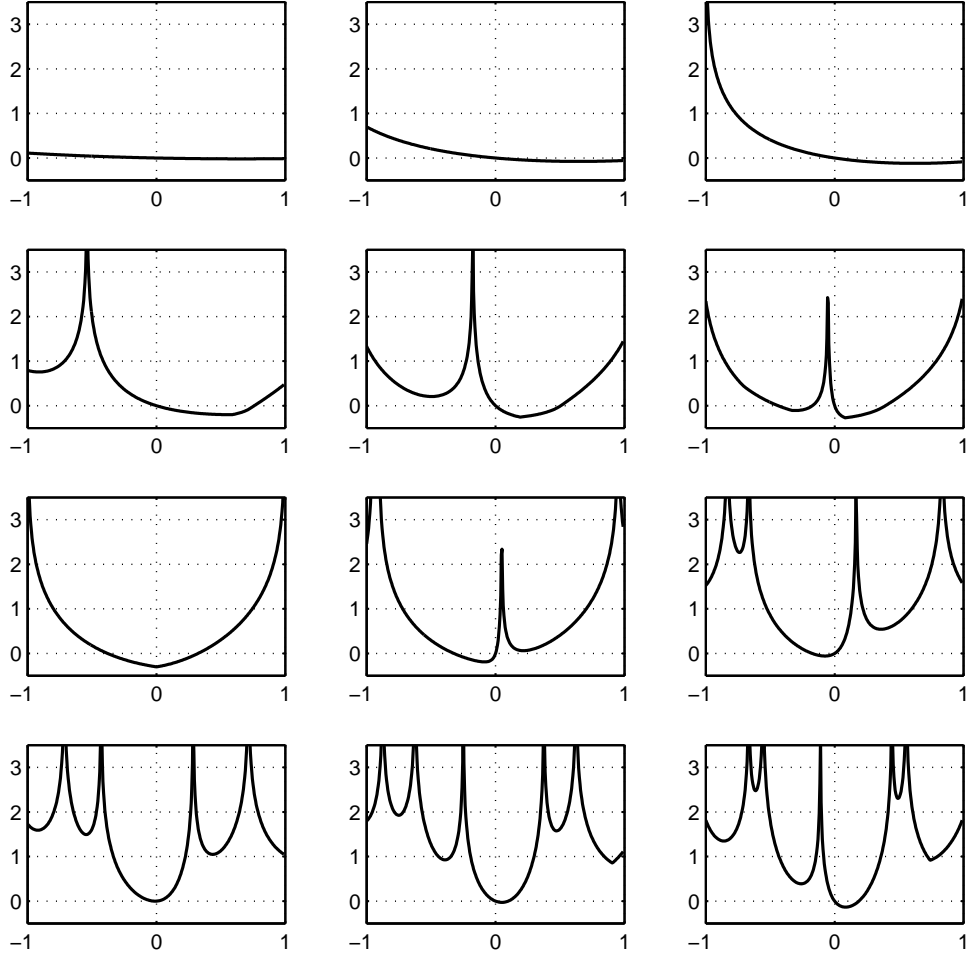}}}
\caption{$\log_{10}(\beta_{\sigma,\alpha})$ vs. $\alpha$ for 
various $\sigma$. From left to right, top to bottom, 
$\sigma/\pi = 
0.2, 0.4, 0.5, 0.65, 0.85, 0.95, 1, 1.05, 1.2, 1.4, 1.6, 1.8$.}
\end{figure} 
%
\section{Preliminary transform in radius}\label{Preliminaries}    

From the pressure formula (\ref{E.p}) we see that the commutator $[\Delta,P] u$
represents the contribution of the viscosity term to the Navier-Stokes pressure
gradient.  Specifically, $[\Delta,P] u$ represents the pressure gradient for
the linear Stokes equations with no-slip boundary and without forcing.  
For this reason, as in \cite{LLP}, we refer to the corresponding pressure 
as the {\em Stokes pressure}, denoted $\ps=\ps(u)$.  
From (\ref{helmholtz}), when $a=u\in H^2(\Omega)$ with $\Omega$ unbounded, we have
$\grad\lap q= \lap\grad q = \grad\div a$ and it follows easily (as in \cite{LLP})
that 
\[
[\Delta,P]u = (I-P)(\Delta u - \grad\grad\cdot u) = \grad \ps.
\]
We recall from \cite[Sec.~2.1]{LLP} that the Stokes pressure $\ps$ is determined
(up to constant) as the solution to the boundary value problem
\begin{equation}\label{BVP}
\Delta \ps =0 
\quad\text{ in } \Omega, 
\qquad
n\cdot \nabla \ps = n\cdot (\Delta-\nabla\nabla\cdot)u 
\quad\text{ on }\Gamma.
\end{equation}
(The boundary condition holds in $H^{-1/2}(\Gamma)$ due to a standard
trace theorem, since the vector fields $\lap u-\grad \div u$ 
and $\grad\ps$ are in $L^2(\Omega,\R^N)$ with zero divergence.)

Letting   
\begin{equation*}
\begin{split}
&I_p=\|\nabla p_s\|^2_{V^0_{2,\alpha}}= 
\int_{\mathcal{K}_{\sigma}}{r^{2\alpha}|\nabla p_s|^2r}\,dr\,d\theta,\\
&I_u= \|\Delta u\|^2_{V^0_{2,\alpha}} 
= \int_{\mathcal{K}_{\sigma}} {r^{2\alpha}|\Delta u|^2r} \, dr\,d\theta,
\end{split}
\end{equation*}
we see that in order to prove Theorem \ref{main}, we must determine
the smallest constant $\beta_{\sigma,\alpha}$ satisfying the inequality
$I_p\leq \beta_{\sigma,\alpha} I_u$, subject to ($\ref{BVP}$).  In this
section, we perform the first steps in our attempt to find $\beta_{\sigma,\alpha}$.
These steps amount to taking a Mellin transform of the problem.
We first rewrite $I_p$, $I_u$ and ($\ref{BVP}$) in terms of the polar
coordinates $(r,\theta)$, then change variables using $r=e^s$, which transforms
$\mathcal{K}_{\sigma}$ to an infinite strip $S$.
Taking a Fourier transform will reduce the problem to a family
of maximization problems parametrized by a Fourier variable $k\in\R$.

\subsection{}
We begin by letting
\[ J=\left( \begin{array}{cc}
0 & -1  \\
1 & 0   \\
\end{array} \right) 
, \qquad
e^{J\theta}=\left( \begin{array}{cc}
\cos \theta & -\sin\theta  \\
\sin\theta & \cos \theta   \\
\end{array}\right).\] 
A straightforward calculation shows that 
\begin{align*}
    \begin{matrix}
        \nabla p_s = r^{-1} e^{J\theta} \left(
           \begin{matrix}
           r\partial_r p\\
           \partial_{\theta} p 
           \end{matrix}   \right),         
   \end{matrix}
\end{align*}
allowing us to rewrite $I_p$ as 
\begin{equation*}
I_p=\int_{\mathcal{K}_{\sigma}} \left( |r\partial_r p|^2 + |\partial_{\theta} p|^2 \right) 
r^{2\alpha-1} \, dr\, d\theta.
\end{equation*}
We change variables by letting $r=e^s$, resulting in a transformation
of the domain $\mathcal{K}_{\sigma}$ to an infinite strip $S=\{(s,\theta)\in\R^2 :
-\infty<s<\infty,$ $0<\theta<\sigma\}$.  We then let $q=e^{\alpha s}p_s$ 
and express $I_p$ in terms of $q$.  We conclude that
\begin{equation}\label{Ip}
I_p = \int_S{\left(|\partial_sq-\alpha q|^2+|\partial_{\theta}q|^2\right)}\, ds\,d\theta = \int_{-\infty}^{\infty} I_{p,k} \, dk,
\end{equation}
where $k$ is the Fourier variable corresponding to $s$, and 
\begin{equation}\label{Ipk1}
I_{p,k} = \int_0^{\sigma} \left( | (k+i\alpha)\hat{q} |^2  +  | \partial_{\theta} \hat{q}|^2 \right) \, d\theta.
\end{equation}

\subsection{}
To rewrite $I_u$, we first calculate
\begin{equation*}
\Delta u = \nabla\cdot\nabla u=
(r\partial_r+2)(r^{-1}\partial_r u) + 
\partial_{\theta}(r^{-2}\partial_{\theta} u)
= r^{-2}( (r\partial_r)^2 +\partial^2_{\theta})u.
\end{equation*}
If we let $u=re^{J\theta}v$, we can show that
\begin{equation}\label{laplaceu}
\Delta u = r^{-1}e^{J\theta} \left( (r\partial_r+1)^2v +(\partial_{\theta}+J)^2v \right).
\end{equation}
We again change variables to express $I_u$ as an integral over $S$.
We let $w=e^{s\alpha}v$, and we find that
\begin{equation}\label{Iu}
\begin{split}
I_u&=\int_S{e^{2s\alpha } |(\partial_s+1)^2 v +(\partial_{\theta}
+J)^2v|^2}\, ds\,d\theta\\
& =\int_S{ |(\partial_s+1-\alpha)^2 w +(\partial_{\theta} +J)^2w|^2}\,
ds\,d\theta
=\int_{-\infty}^{\infty}  I_{u,k} \, dk,
\end{split}
\end{equation}
where
\begin{equation}\label{Iuk1}
I_{u,k} = \int_0^{\sigma} | (ik+1-\alpha)^2 \hat{w} + (\partial_{\theta} + J)^2\hat{w} |^2 \,d\theta.
\end{equation}

\subsection{}
As with $I_p$ and $I_u$, we wish to rewrite ($\ref{BVP}$) in terms of
$k$, $\theta$, $q$, and $w$.  We perform a change of variables and
rewrite the first condition of ($\ref{BVP}$) as
$e^{-2s}(\partial_s^2+\partial_{\theta}^2)p=0$.  Recalling that
$q=e^{s\alpha}p$, we find that $\Delta q-e^{-2s}(2\alpha\partial_s
-\alpha^2 )q=0$, so
$(\partial^2_{\theta}+\partial^2_s-2\alpha\partial_s  +\alpha^2) q=0$.
To rewrite the boundary condition of ($\ref{BVP}$), we observe that
the left hand side can be rewritten using the equalities $\langle
e_{\theta} ,\nabla p_s\rangle=\langle e_2, r^{-1}(r\partial_r
p_s,\partial_{\theta}p_s) \rangle= r^{-1}\partial_{\theta} p_s$.  For
the right hand side, we use ($\ref{laplaceu}$), combined with the
equality
\begin{align*}
    \begin{matrix}
         \nabla\nabla \cdot u=r^{-1}e^{J\theta} \left(
           \begin{matrix}
           r\partial_r \nabla\cdot u\\
           \partial_{\theta}\nabla\cdot u 
           \end{matrix}   \right) = r^{-1}e^{J\theta}\left(\begin{matrix}  
           r\partial_r((r\partial_r +2)v_1 +\partial_{\theta}v_2)\\
           \partial_{\theta}((r\partial_r +2)v_1 +\partial_{\theta}v_2)
              \end{matrix}         \right)
   \end{matrix}
\end{align*}
and the property $v=0$ when $\theta=0$ and $\theta=\sigma$, to conclude that 
\begin{equation*}
\begin{split}
&\langle e_{\theta}, \Delta u-\nabla\nabla\cdot u\rangle
=r^{-1}((r\partial_r)^2 +2r\partial_r)v_2-\partial_{\theta}\partial_r v_1\\
&\qquad\qquad= -r^{-2}v_2-\partial_{\theta}\partial_r v_1 
= -\partial_{\theta}\partial_r v_1
\end{split}
\end{equation*} 
for $\theta=0$ and $\theta=\sigma$.  We can therefore rewrite the
boundary condition in ($\ref{BVP}$) as 
\begin{equation}
\partial_{\theta} p_s=-\partial_s\partial_{\theta}v_1.  
\end{equation}
Using the equality
$w=e^{s\alpha}v$, we see after a calculation that we can recast
($\ref{BVP}$) as the following boundary value problem on $S$:
\begin{equation}\label{BVP1}
\begin{split}
&(\partial^2_{\theta}+\partial^2_s-2\alpha\partial_s  +\alpha^2) q
=0 \quad\text{ in } S, \\
&\partial_{\theta} q 
= -\partial_{\theta}(\partial_s-\alpha)w_1\quad\text{ when } \theta=0,\sigma.
\end{split}
\end{equation}
Finally, taking the Fourier transform of ($\ref{BVP1}$) in $s$, we have 
that for each $k\in\R$, $\hat q$ must solve the boundary value problem
\begin{equation}\label{harmonic}
\begin{split}
&\partial^2_{\theta} \hat{q} = (k+i\alpha)^2 \hat{q} \quad\text{ for } 0<\theta<\sigma,\\
&\partial_{\theta} \hat{q} = -\partial_{\theta} (ik-\alpha){\hat{w}}_1
\quad\text{ when } \theta=0,\sigma.
\end{split}
\end{equation}

\section{Optimization in angle}\label{solveforbeta}  
In this section, we determine $\beta_{\sigma,\alpha}=\sup
\frac{I_p}{I_u}$ subject to ($\ref{harmonic}$) and the no-slip
boundary condition.
First, for $k\ne0$ we suppress the $\alpha$ and $\sigma$ variables and define
\begin{equation}\label{bdefs}
\beta_{k}=\sup\left\{ \frac{I_{p,k}}{I_{u,k}}:
\text{(\ref{harmonic}) holds, and $\hat w=0$ for $\theta=0,\sigma$}\right\}. 
\end{equation}
Note that since $w$ is real, we have 
$\hat w(-k,\theta)=\overline{\hat w(k,\theta)}$,
hence $I_{u,-k}=I_{u,k}$ from (\ref{Iuk1}), 
and similarly $I_{p,-k}=I_{p,k}$ from (\ref{Ipk1}). 
We conclude that $\beta_k$ is even in $k$.  We define
${\hat{\beta}}_{\sigma,\alpha} =\sup_{k>0} {{\beta}_k}$, and we observe that 
\begin{equation}\label{betasigmak}
I_p=\int_{-\infty}^{\infty} {I_{p,k}}\, dk  
\leq  \int_{-\infty}^{\infty} {\beta_{k} I_{u,k}} \,dk 
\leq {\hat{\beta}}_{\sigma,\alpha} I_u.
\end{equation}
We will prove Theorem 3 by computing that $\beta_{k}=\max\{\hat\beta_{+,k},\hat\beta_{-,k}\}$ 
as given by (\ref{E.bpmk}), 
and by showing that 
$\hat\beta_{\sigma,\alpha} \le {\beta}_{\sigma, \alpha}$.  
Since evidently 
$\hat\beta_{\sigma,\alpha} \ge {\beta}_{\sigma, \alpha}$,
the result will follow.  

\subsection{}
We first rewrite the quantity $I_{u,k}$ from ($\ref{Iuk1}$) to diagonalize the matrix
involved.  We define
\[ V=\left( \begin{array}{cc}
1 & 1  \\
i & -i   \\
\end{array} \right) ,\qquad
 \Lambda=\left( \begin{array}{cc}
-1 & 0  \\
0 & 1   \\
\end{array}\right).\] 
Then letting $-i\hat{w}=Vy$ with $y=(y_1,y_2)$, and using 
$JV=V(i\Lambda)$, we rewrite $I_{u,k}$ in the following way: 
\begin{eqnarray}
I_{u,k} &=& \int_0^\sigma{ |\left((ik+1-\alpha)^2  +(\partial^2_{\theta}
+2J\partial_{\theta}-1)\right)(Vy)|^2}\,d\theta
\nonumber\\
&=& 2\int_0^{\sigma}{ |\left((ik+1-\alpha)^2  +(\partial^2_{\theta}
+2i\Lambda\partial_{\theta}-1)\right)y|^2}\,d\theta
\nonumber\\
&=& 
2\int_0^{\sigma}{(|L_1y_1|^2+|L_2y_2|^2) }\,d\theta,
\label{iuk}
\end{eqnarray}
where 
\begin{eqnarray*}
L_1&=&(ik+1-\alpha)^2+\partial^2_{\theta}-2i\partial_{\theta}-1, \\
L_2&=&(ik+1-\alpha)^2+\partial^2_{\theta}+2i\partial_{\theta}-1 .
\end{eqnarray*}

\subsection{}
We next express the quantity $I_{p,k}$ from ($\ref{Ipk1}$)
in terms of the boundary data from (\ref{harmonic}).  
From ($\ref{harmonic}$) it is clear that explicitly
\begin{equation}\label{qhat}
\begin{split}
&\hat{q}(k,\theta)
=\alpha_+e^{(k+i\alpha)(\theta-\sigma)}+\alpha_-e^{-(k+i\alpha)\theta},\\
&\partial_{\theta} \hat{q}(k,\theta) 
= (k+i\alpha)\alpha_+e^{(k+i\alpha)(\theta-\sigma)} 
- (k+i\alpha)\alpha_-e^{-(k+i\alpha)\theta},
\end{split}
\end{equation} 
for some complex constants $\alpha_+$ and $\alpha_-$.  If we define
\begin{equation}\label{D.om}
\omega=e^{-(k+i\alpha)\sigma}
\end{equation}
for convenience, we see from ($\ref{qhat}$) that
 \begin{equation}\label{qhatboundary}
    \begin{matrix}
       \hat{q}(k,\theta)= \left\{ 
       \begin{matrix}
        \alpha_+ + \alpha_-\omega, & \theta=\sigma,\\ 
        \alpha_+\omega+\alpha_-, &\theta=0,\\
        \end{matrix}  
        \right.
 \end{matrix}
 \end{equation}
and 
\begin{equation}\label{partialqhat}
    \begin{matrix}
       \partial_{\theta} \hat{q}(k,\theta)= \left\{
            \begin{matrix}
                (k+i\alpha)(\alpha_+-\alpha_-\omega),  & \theta=\sigma, \\
                 (k+i\alpha)(\alpha_+\omega-\alpha_-),  &\theta=0. \\
            \end{matrix}
            \right.
    \end{matrix}
\end{equation}
Combining (\ref{partialqhat}) with the equality
$-i{\hat{w}}_1=(Vy)_1$, we can rewrite the boundary conditions in
($\ref{harmonic}$) as 
\begin{equation}\label{yboundary}
    \begin{matrix}
       \alpha_+-\alpha_-\omega=\partial_{\theta}(y_1+y_2),  & \theta=\sigma, \\
  \alpha_+\omega-\alpha_- =\partial_{\theta}(y_1+y_2),  &\theta=0. \\
            \end{matrix}
\end{equation}
These equations will be used later to determine $\alpha_+$ and $\alpha_-$ 
from $y$
(note $\omega^2\ne1$).
To rewrite $I_{p,k}$, we apply (\ref{harmonic}) and integrate by parts.  This
gives
\\
\begin{equation*}
\begin{split}
&\int_S{|(\partial_s-\alpha)q|^2 }\,ds\,d\theta=\int_S
{(k+i\alpha)\hat{q}(k-i\alpha)\bar{\hat{q}}}\,dk\,d\theta\\
&= \int_{-\infty}^{\infty}{\left(
\partial_{\theta}\hat{q}(\sigma) \bar{\hat{q}}(\sigma) 
-\partial_{\theta}\hat{q}(0)\bar{\hat{q}}(0)\right)}\,dk - \int_S \partial_{\theta}\hat{q}\partial_{\theta}\bar{\hat{q}}\,dk\,d\theta\\
&\qquad\qquad -\int_S 2(ik\alpha-\alpha^2)\hat{q}\bar{\hat{q}}\,dk\,d\theta,
\end{split}
\end{equation*}
which, in light of ($\ref{Ip}$), allows us to write 
\begin{equation}\label{findIpk}
I_{p,k}=
\partial_{\theta}\hat{q}(\sigma)\bar{\hat{q}}(\sigma) 
-\partial_{\theta}\hat{q}(0)\bar{\hat{q}}(0)
-\int_0^{\sigma} 2(ik\alpha-\alpha^2)\hat{q}\bar{\hat{q}}\,d\theta.
\end{equation} 
In order to write $\int_0^{\sigma}
2(ik\alpha-\alpha^2)\hat{q}\bar{\hat{q}}\,d\theta$ in terms of
$\alpha_+$ and $\alpha_-$, we use ($\ref{qhat}$) 
to evaluate the dot product and integrate.  
We conclude that 
\begin{equation}\label{Ip1}
\begin{split}
&\int_0^{\sigma} 2(ik\alpha-\alpha^2)\hat{q}\bar{\hat{q}}\,d\theta 
= (ik\alpha-\alpha^2) 
\left( \frac{|\alpha_+|^2}{k}+\frac{|\alpha_-|^2}{k}\right)
\left(1-e^{-2k\sigma}\right)\\
&\qquad \qquad 
+ (ik-\alpha)
\left(i\alpha_-{\bar{\alpha}}_+ + i\alpha_+{\bar{\alpha}}_- \right)
e^{-(k+i\alpha)\sigma}\left(1-e^{2i\alpha\sigma}\right)\\
& \qquad\qquad= (ik\alpha-\alpha^2) 
\left(\frac{|\alpha_+|^2}{k}+\frac{|\alpha_-|^2}{k}\right)
\left(1-|\omega|^2\right)\\
&\qquad\qquad + (k+i\alpha)
(\alpha_-{\bar{\alpha}}_+ + \alpha_+{\bar{\alpha}}_-)
(\bar{\omega} - \omega).
\end{split}
\end{equation} 
Similarly, to compute
$\partial_{\theta}\hat{q}(\sigma)\bar{\hat{q}}(\sigma)
-\partial_{\theta}\hat{q}(0)\bar{\hat{q}}(0)$, we use the
formulas for $\hat{q}$ and $\partial_{\theta}\hat{q}$ on the boundary
given in ($\ref{qhatboundary}$) and ($\ref{partialqhat}$) to write
\begin{equation}\label{Ip2}
\begin{split}
&\partial_{\theta}\hat{q}(\sigma)\bar{\hat{q}}(\sigma) 
-\partial_{\theta}\hat{q}(0)\bar{\hat{q}}(0)
= (k+i\alpha) \{ (|\alpha_+|^2+|\alpha_-|^2)(1-|\omega|^2)\\
&\qquad\qquad\qquad
+(\alpha_+{\bar{\alpha}}_- +\alpha_-{\bar{\alpha}}_+ )
(\bar{\omega}-\omega) \}. 
\end{split}
\end{equation}
Plugging ($\ref{Ip1}$) and ($\ref{Ip2}$) into ($\ref{findIpk}$), we
discover that 
\begin{equation*}
I_{p,k} = \left(   \frac{k^2+\alpha^2}{k}  \right)  
\left( |\alpha_+|^2+|\alpha_-|^2 \right)   
\left(1-|\omega|^2 \right)  .
\end{equation*}

\subsection{}
By the results of the previous subsection, $\beta_k$ is the supremum of 
the ratio $I_{p,k}/I_{u,k}$ subject to (\ref{yboundary}) 
and the no-slip boundary conditions $y=0$ for $\theta=0,\sigma$.
In order to compute $\beta_k$, we will argue that the supremum in 
(\ref{bdefs}) is a maximum and use a variational argument. 

The existence of a maximizer is proved by a standard argument 
in the calculus of variations:
It is clear that $0<\beta_k\le\infty$, and that the 
ratio $I_{p,k}/I_{u,k}$ is a homogeneous function of $y$. 
Thus we may choose a maximizing sequence of vector functions $y$ 
with fixed $H^2$ Sobolev norm on $[0,\sigma]$.
Evidently, the quantities $L_1y_1$, $L_2y_2$ remain bounded in 
$L^2$, and the complex scalar quantities $\partial_\theta y_1$,
$\partial_\theta y_2$ at $\theta=0,\sigma$ remain bounded.
We may choose a subsequence converging weakly in $H^2$ 
such that the quantities
$\partial_\theta y_1$, $\partial_\theta y_2$ 
at $\theta=0,\sigma$ converge.
Then the weak limit is a maximizer by weak lower semicontinuity of
the $L^2$ norm.

Next, consider any smooth curve $\tau\mapsto y=y(\tau)$ into $H^2$ with the
property that (\ref{yboundary}) and the no-slip conditions hold for all $\tau$,
and $I_{p,k}/I_{u,k}$ achieves its maximum at $\tau=0$.
Then at $\tau=0$ we have
\begin{equation}\label{E.maxb}
0={\dot{I}}_{p,k}-\beta_{k}{\dot{I}}_{u,k}.
\end{equation}  
We now determine ${\dot{I}}_{p,k}$ and ${\dot{I}}_{u,k}$ and solve for
$\beta_{k}$.  Differentiating $I_{p,k}$, we find
\begin{equation*}
{\dot{I}}_{p,k} = \left( \frac{k^2+\alpha^2}{k} \right)
({\bar{\dot{\alpha}}}_+{\alpha}_+ + {\bar{\alpha}}_+ {\dot{\alpha}}_++
{\bar{\dot{\alpha}}}_-{\alpha}_- + {\bar{\alpha}}_-
{\dot{\alpha}}_-)(1-|\omega|^2).
\end{equation*}
From ($\ref{yboundary}$) we infer that
\begin{equation}\label{solvefor+and-}
\begin{split}
&\alpha_+(1-\omega^2)=\partial_{\theta}
(y_1+y_2)e^{(k+i\alpha)(\theta-\sigma)}|^{\sigma}_0,\\
&\alpha_-(1-\omega^2)=\partial_{\theta}
(y_1+y_2)e^{-(k+i\alpha)\theta}|^{\sigma}_0.
\end{split}
\end{equation}
By differentiating in $\tau$, we can solve for ${{\dot{\alpha}}}_+$ and
${{\dot{\alpha}}}_-$, allowing us to eliminate
${{\dot{\alpha}}}_+$ and ${{\dot{\alpha}}}_-$ from the formula
for ${\dot{I}}_{p,k}$.  Indeed, if we let
$\gamma_i(\theta)=\partial_{\theta}{\dot{y}}_i$ for $i=1$, $2$, 
we have
\begin{equation*}
\begin{split}
&\bar{\dot{\alpha}}_+=
\frac{(1-\omega^2)}{|1-\omega^2|^2}({\bar{\gamma}}_1
+{\bar{\gamma}}_2)e^{(k-i\alpha)(\theta-\sigma)}|^{\sigma}_0,\\
&\bar{\dot{\alpha}}_-=\frac{(1-\omega^2)}{|1-\omega^2|^2}({\bar{\gamma}}_1
+{\bar{\gamma}}_2)e^{-(k-i\alpha)\theta}|^{\sigma}_0.
\end{split}
\end{equation*}

Similarly, we differentiate $I_{u,k}$. 
Letting $L_1^*$ and $L_2^*$ denote the formal adjoints of
$L_1$ and $L_2$, respectively, and recalling from the no-slip boundary
conditions that $\dot{y}=0$ at $\theta=0,\sigma$, 
we integrate by parts to conclude that 
\begin{equation*}
{\dot{I}}_{u,k} = 4\left(\int_0^{\sigma}{\Re(\rho)} \, d\theta 
+ \Re(\psi) \right),
\end{equation*}
where 
\begin{equation*}
\rho=\bar{{\dot{y}}}_1 (L_1^*L_1)y_1 + \bar{{\dot{y}}}_2(L_2^*L_2)y_2,
\qquad \psi=\overline{\partial_{\theta}{\dot{y}}}_1(L_1y_1) 
+ \overline{\partial_{\theta}{\dot{y}}}_2(L_2y_2) |^{\sigma}_0.
\end{equation*}
We observe that 
\begin{equation*}
L_1L_1^* y_1=0 \quad\text{ and }\quad L_2L_2^* y_2=0, 
\end{equation*}
so that
${\dot{I}}_{u,k}$ reduces to ${\dot{I}}_{u,k}=4\Re( \psi )$.

Using this information, we
can rewrite (\ref{E.maxb}) as
\begin{equation}\label{23}
\begin{split}
0 &= 2\Re\left\{ {\bar{\gamma}}_1 \left( 2\beta_k L_1y_1 - 
\left( \frac{k^2+\alpha^2}{k} \right)
\frac{(1-|\omega|^2)}{1-{\bar{\omega}}^2}
(\alpha_+e^{(k-i\alpha)(\theta-\sigma)}
+ {\alpha}_-e^{-(k-i\alpha)\theta} )\right)|^{\sigma}_0\right\}\\
& +2\Re\left\{ {\bar{\gamma}}_2 \left( 2\beta_k L_2y_2 - 
\left( \frac{k^2+\alpha^2}{k} \right)
\frac{(1-|\omega|^2)}{1-{\bar{\omega}}^2}
(\alpha_+e^{(k-i\alpha)(\theta-\sigma)}+
{\alpha}_-e^{-(k-i\alpha)\theta} )\right)|^{\sigma}_0\right\}.     
\end{split}
\end{equation} 
Since $\gamma_1(\theta)$ and $\gamma_2(\theta)$ are arbitrary at
$\theta=0$ and $\theta=\sigma$, ($\ref{23}$) yields four (natural) boundary
conditions:
\begin{equation}\label{27}
\begin{split}
&2\beta_k L_1y_1 = \left( \frac{k^2+\alpha^2}{k} \right)
\frac{(1-|\omega|^2)}{1-{\bar{\omega}}^2}
(\alpha_+ + {\alpha}_-\bar{\omega} ) 
\text{ and }\\
&2\beta_k L_2y_2 = \left( \frac{k^2+\alpha^2}{k} \right)
\frac{(1-|\omega|^2)}{1-{\bar{\omega}}^2}
(\alpha_+ + {\alpha}_-\bar{\omega}  ), \text{ when } \theta=\sigma,\\
&2\beta_k L_1y_1 = \left( \frac{k^2+\alpha^2}{k} \right)
\frac{(1-|\omega|^2)}{1-{\bar{\omega}}^2}
(\alpha_+\bar{\omega} + {\alpha}_- ) \text { and }\\
&2\beta_k L_2y_2 = \left( \frac{k^2+\alpha^2}{k} \right)
\frac{(1-|\omega|^2)}{1-{\bar{\omega}}^2}
(\alpha_+\bar{\omega} + {\alpha}_- ), \text{ when } \theta=0.
\end{split}
\end{equation}
In addition, we have the four no-slip boundary conditions
\begin{equation}\label{BC1}
y_1(\sigma)=y_2(\sigma)=y_1(0)=y_2(0)=0.
\end{equation}

\subsection{}
Using ($\ref{BC1}$), ($\ref{27}$), ($\ref{yboundary}$), and the
property $L_1^*L_1y_1=L_2^*L_2y_2=0$ on $(0,\sigma)$, we can explicitly 
solve for the maximizer of $\beta_k$. 
To simplify the calculations in what follows, we first use reflection
symmetry to show that either 
\begin{equation*}
(\alpha_+,\alpha_-)=(1,1) \quad\text{ or }\quad (\alpha_+,\alpha_-)=(1,-1).
\end{equation*}  
Letting
$\hat{\theta}=\sigma-\theta$, we see from our construction of
$\hat{q}$ in ($\ref{qhat}$) that $\alpha_+$ and $\alpha_-$ exchange
roles after reflection; thus, it is natural to set
${\hat{\alpha}}_+=\alpha_-$, and ${\hat{\alpha}}_-=\alpha_+$.   In
addition, we let ${\hat{y}}_2(\hat{\theta})=y_1(\theta)$ and
${\hat{y}}_1(\hat{\theta})=y_2(\theta)$.  A straightforward
calculation shows that $({\hat{y}}_1, {\hat{y}}_2, {\hat{\alpha}}_+,
{\hat{\alpha}}_-)$ solves the set of linear equations consisting of
($\ref{BC1}$), ($\ref{27}$), ($\ref{yboundary}$), and
$L_1^*L_1y_1=L_2^*L_2y_2=0$.\Obsolete{Specifically, 
\begin{equation*}
\begin{split}
&\partial_{\theta}({\hat{y}}_1 + {\hat{y}}_2)(0)=  
-\alpha_+ +\alpha_-\omega = {\hat{\alpha}}_+\omega -{\hat{\alpha}}_-, 
\text{ and}\\
&\partial_{\theta}({\hat{y}}_1 + {\hat{y}}_2)(\sigma)=
-\alpha_+\omega + \alpha_- = {\hat{\alpha}}_+ -{\hat{\alpha}}_-\omega .
\end{split}
\end{equation*}
Moreover,
\begin{equation*}
\begin{split}
&2\beta L_1{\hat{y}}_1 = \left( \frac{k^2+\alpha^2}{k} \right)
\frac{(1-|\omega|^2)}{1-{\bar{\omega}}^2}
({\hat{\alpha}}_+ + {\hat{\alpha}}_-\bar{\omega} ) 
\text{ and }\\
&2\beta L_2{\hat{y}}_2= \left( \frac{k^2+\alpha^2}{k} \right)
\frac{(1-|\omega|^2)}{1-{\bar{\omega}}^2}
({\hat{\alpha}}_+ + {\hat{\alpha}}_-\bar{\omega}  ), 
\text{ when } \hat{\theta}=\sigma.\\
&2\beta L_1{\hat{y}}_1 = \left( \frac{k^2+\alpha^2}{k} \right)
\frac{(1-|\omega|^2)}{1-{\bar{\omega}}^2}
({\hat{\alpha}}_+\bar{\omega} + {\hat{\alpha}}_- ) 
\text { and }\\
&2\beta L_2{\hat{y}}_2 = \left( \frac{k^2+\alpha^2}{k} \right)
\frac{(1-|\omega|^2)}{1-{\bar{\omega}}^2}
({\hat{\alpha}}_+\bar{\omega} + {\hat{\alpha}}_- ), 
\text{ when } \hat{\theta}=0.
\end{split}
\end{equation*}}
We deduce that 
\[(y_1+{\hat{y}}_1, y_2+{\hat{y}}_2, \alpha_+ +
\alpha_-, \alpha_- +\alpha_+) 
\quad\text{ and }\quad
(y_1-{\hat{y}}_1, y_2-{\hat{y}}_2, \alpha_+-\alpha_-, \alpha_- -\alpha_+)
\] also solve these equations.
We conclude that every pair $(\alpha_+, \alpha_-)$ will yield the same
value for $\beta_k$ as either $(\alpha_+,\alpha_-)=(1,1)$ or
$(\alpha_+,\alpha_-)=(1,-1)$.  Therefore it suffices to consider only these cases.  

\subsection{}
We can eliminate $y_2$ by
observing that if $(\alpha_+,\alpha_-)=(1,1)$, then
$y_2(\theta)={\hat{y}}_2(\theta)=y_1(\sigma-\theta)$, and if
$(\alpha_+,\alpha_-)=(1,-1)$, then
$y_2(\theta)=-{\hat{y}}_2(\theta)=-y_1(\sigma-\theta)$. 
Then we infer
from boundary conditions in ($\ref{yboundary}$) that
\begin{equation}\label{reducetoy1}
\begin{split}
&1-\omega=\partial_{\theta} y_1(\sigma) - \partial_{\theta} y_1(0), 
\text{ when }(\alpha_+,\alpha_-)=(1,1), 
\text{ and}\\
&1+\omega=\partial_{\theta} y_1(\sigma) + \partial_{\theta} y_1(0), 
\text{ when }(\alpha_+,\alpha_-)=(1,-1).
\end{split}
\end{equation}  
We are now in a position to solve for $y_1$ and ultimately $\beta_k$.
We first recall that
$L_1=(ik+1-\alpha)^2+\partial^2_{\theta}-2i\partial_{\theta}-1$\Obsolete{and
$L_2=(ik+1-\alpha)^2+\partial^2_{\theta}+2i\partial_{\theta}-1$},
while, formally, the adjoint of this operator is given by
$L_1^*=(ik-1+\alpha)^2+\partial^2_{\theta}-2i\partial_{\theta}-1$\Obsolete{
and $L^*_2=(ik-1+\alpha)^2+\partial^2_{\theta}+2i\partial_{\theta}-1$,
respectively}.  The characteristic polynomials of these two operators
are
\begin{equation*}
\begin{split}
&p_1(\mu) = (\mu - (2i-k-i\alpha))(\mu-(k+i\alpha)),\\
&p^{\ast}_1(\mu) = (\mu-(2i+k-i\alpha))(\mu-(-k+i\alpha)).
\end{split}
\end{equation*}   
Since $L^{\ast}_1L_1y_1=0$ on $(0,\sigma)$, 
we can conclude that $y_1(\theta)$ takes the form
\begin{equation}\label{26}
\begin{split}
y_1(\theta) = a_1e^{(k+i\alpha)(\theta-\sigma)} +
a_2e^{-(k-2i+i\alpha)\theta} + a_3e^{-(k-i\alpha)\theta}+
a_4e^{(k+2i-i\alpha)(\theta-\sigma)}\\
\end{split}
\end{equation}
for some constants $a_i$, $1\leq i \leq 4$.  The boundary conditions
$y_1(\sigma)=y_1(0)=0$, combined with ($\ref{26}$), yield the two
equalities
\\
\begin{equation}\label{29}
\begin{split}
&0=a_1+a_2\omega e^{2i\sigma} + a_3\bar\omega + a_4,\\
&0=a_1\omega + a_2 + a_3 + a_4\bar\omega e^{-2i\sigma}.
\end{split}
\end{equation}
We will use the boundary conditions for $y_1$ in ($\ref{27}$) combined
with the equalities in ($\ref{29}$) to write the four unknowns $a_j$,
$1\leq j \leq 4$, in terms of $\alpha_+$ and $\alpha_-$.  

Using the equality $L_j=L^{\ast}_j + 4(1-\alpha)ik$ for $j=1$,
$2$, and ($\ref{26}$), we conclude that
\[
L_1y_1=4(1-\alpha)ik\left(a_3e^{-(k-i\alpha)\theta} +
a_4e^{(k+2i-i\alpha)(\theta-\sigma)}\right) .
\]
Plugging this information into the
two boundary conditions in ($\ref{27}$) yields the two equalities
\begin{equation}\label{28}
\begin{split}
&8\beta_k (1-\alpha)ik(a_3\bar\omega +a_4) = \left(
\frac{k^2+\alpha^2}{k}
\right)\frac{(1-|\omega|^2)}{1-{\bar{\omega}}^2}(\alpha_+ +
\alpha_-\bar{\omega}),\\
&8\beta_k(1-\alpha)ik(a_3 + a_4\bar\omega e^{-2i\sigma}) =
\left( \frac{k^2+\alpha^2}{k}
\right)\frac{(1-|\omega|^2)}{1-{\bar{\omega}}^2}(\alpha_-+\alpha_+\bar{\omega}
).
\end{split}
\end{equation}

\subsection{}
The value of $\beta_k$ is determined by the equations in
(\ref{28}) and (\ref{29}) together with (\ref{reducetoy1}).
Evidently $\beta_k=\max\{\bkp,\bkm\}$ where $\bkp$ and $\bkm$ 
are the values determined from these equations in 
each of the two cases
$(\alpha_+,\alpha_-)=(1,1)$ and $(\alpha_+,\alpha_-)=(1,-1)$
respectively.  

With $(\alpha_+,\alpha_-)=(1,1)$, using the four equations given in
(\ref{28}) and (\ref{29}), we solve for the unknowns
$a_j$, $1\leq j \leq 4$, finding that
\begin{equation}\label{avalue}
\begin{split}
&a_1\bkp(1-\omega^2 e^{2i\sigma})=-\phi_1(1-\omega e^{2i\sigma}),\\
&a_2\bkp(1-\omega^2 e^{2i\sigma})=-\phi_1(1-\omega),\\
&a_3\bkp(1-{\bar{\omega}}^2e^{-2i\sigma})=\phi_1(1-\bar{\omega}e^{-2i\sigma}),\\
&a_4\bkp(1-{\bar{\omega}}^2e^{-2i\sigma})=\phi_1(1-\bar{\omega}),
\end{split}
\end{equation}
where
\begin{equation*}
\phi_1=\frac{(k^2+\alpha^2)(1-|\omega|^2)(1+\bar{\omega})}
{8i k^2(1-\alpha)(1-{\bar{\omega}}^2)}.
\end{equation*}
Using (\ref{reducetoy1}) with (\ref{26}), we see that
\begin{equation}\label{41}
\begin{split}
&1-\omega=a_1\hat{k}(1-\omega) 
+ a_2(2i-\hat{k})(\omega e^{2i\sigma}-1)\\
&\qquad\qquad 
+ a_3(\bar{\hat{k}})(1-\bar{\omega})
+ a_4(2i+\bar{\hat{k}})(1-\bar{\omega}e^{-2i\sigma}),
\end{split}
\end{equation}
where $\hat{k}=k+i\alpha$.  Plugging
the formulas for $a_j$ into (\ref{41}) and solving
for $\bkp$ yields 
\begin{equation*}
\bkp =\frac{\phi_1}{(1-\omega)}\left\{ \frac{(1-\omega)(1-\omega
e^{2i\sigma})}{(1-\omega^2e^{2i\sigma})} (2i-2\hat{k}) +
\frac{(1-\bar{\omega})(1-\bar{\omega}
e^{-2i\sigma})}{(1-{\bar{\omega}}^2 e^{-2i\sigma})}
(2i+2\bar{\hat{k}})\right\}.
\end{equation*}  

To find $\bkm$, we let $(\alpha_+,\alpha_-)=(1,-1)$, and we again
use (\ref{28}) and (\ref{29}) to solve for $a_j$,
$1\leq j \leq 4$.  To simplify notation, we define 
\begin{equation*}
\phi_2=\frac{(k^2+\alpha^2)(1-|\omega|^2)(1-\bar{\omega})}{8i
k^2(1-\alpha)(1-{\bar{\omega}}^2)}.
\end{equation*}
We compute the $a_j$ and conclude that
\begin{equation}\label{avalue2}
\begin{split}
&a_1\bkm(1-\omega^2 e^{2i\sigma})=-\phi_2(1+\omega e^{2i\sigma}),\\
&a_2\bkm(1-\omega^2 e^{2i\sigma})=\phi_2(1+\omega),\\
&a_3\bkm(1-{\bar{\omega}}^2e^{-2i\sigma})
=-\phi_2(1+\bar{\omega}e^{-2i\sigma}),\\
&a_4\bkm(1-{\bar{\omega}}^2e^{-2i\sigma})=\phi_2(1+\bar{\omega}).
\end{split}
\end{equation}
We solve for $\bkm$ using (\ref{reducetoy1}) with (\ref{26}) like
before, and find that
\begin{equation*}
\bkm=\frac{\phi_2}{(1+\omega)}\left\{ \frac{(1+\omega)(1+\omega
e^{2i\sigma})}{(1-\omega^2e^{2i\sigma})} (2i-2\hat{k})+
\frac{(1+\bar{\omega})(1+\bar{\omega}
e^{-2i\sigma})}{(1-{\bar{\omega}}^2 e^{-2i\sigma})}
(2i+2\bar{\hat{k}})\right\}.
\end{equation*}
At this point, one can check that $\beta_{\pm,k} = \hat\beta_{\pm,k}$ 
as given in (\ref{E.bpmk}). 

\subsection{}
To complete the proof of Theorem \ref{main}, as indicated at the
beginning of this section, we must show that
$\beta_{\sigma,\alpha}\ge{\hat{\beta}}_{\sigma, \alpha}$.  
To prove this, suppose $\hat\beta<\hat\beta_{\sigma,\alpha}$. 
Then there exists $k_0\ne0$ such that $\beta_{k_0}>\hat\beta$.  
We choose $y$ to be a maximizer of the ratio $I_{p,k_0} / I_{u,k_0}$. 
In a change of notation, we let $I_{p,k}$ and $I_{u,k}$ denote the integrals
corresponding to this fixed $y$, with $\hat{q}$ determined by (\ref{harmonic})
for $k$ varying.  Since $y$ may not be a maximizer for $k\ne k_0$,
we only have $\beta_k\ge I_{p,k}/I_{u,k}$ in general. However, 
by continuity it is evident that there exists $\delta>0$ such that
whenever $|k-k_0|<\delta$ we have $I_{p,k}/I_{u,k}>\hat\beta$. 

Next, we define ${\chi}_{\delta}(k)$ to be a smooth bump function
independent of $\theta$ and supported in a $\delta$-neighborhood
of $k_0$.  Recalling that $-i\hat{w}=Vy$, we set ${\hat{w}}_{\delta}={\chi}_{\delta}\hat{w}$
and ${\hat{q}}_{\delta}={\chi}_{\delta}\hat{q}$, and we observe that
(${\hat{w}}_{\delta},{\hat{q}}_{\delta}$) solves
($\ref{harmonic}$) and ${\hat{w}}_{\delta}=0$ for $\theta=0,\sigma$.  Moreover, if $I_{p_{\delta},k}$ and $I_{u_{\delta},k}$ are the integrals corresponding to ${\hat{w}}_{\delta}$ and ${\hat{q}}_{\delta}$, then one sees that
$I_{p_{\delta},k}={\chi_{\delta}}^2 I_{p,k}$ and
$I_{u_{\delta},k}={\chi_{\delta}}^2 I_{u,k}$.  We can
then write 
\begin{equation}\label{optbeta2}
\frac{I_{p_{\delta}}}{I_{u_{\delta}}} 
= \frac{\int_{k_0-\delta}^{k_0+\delta} 
{\chi_{\delta}}^2 I_{p,k}\,dk }
{ \int_{k_0-\delta}^{k_0+\delta} 
{\chi_{\delta}}^2 I_{u,k}\,dk } 
> \frac{\int_{k_0-\delta}^{k_0+\delta} 
\hat\beta{\chi_{\delta}}^2 I_{u,k}\,dk }
{ \int_{k_0-\delta}^{k_0+\delta} 
{\chi_{\delta}}^2 I_{u,k} \,dk } 
= \hat\beta.
\end{equation}    
We conclude that
$\beta_{\sigma,\alpha}\geq\hat\beta$, hence 
$\beta_{\sigma,\alpha}\geq {\hat{\beta}}_{\sigma,\alpha}$.  

\section{Causes for blowup of the optimal constant}\label{blowupsect}
One can rewrite the formulas for $\bpk$ and $\bmk$ from Theorem \ref{main} in the following way:
\begin{equation}\label{discat0}
\bk{\pm,k}=  \frac{ \psi_1+\psi_2 }{ 2k^2(\cosh (k\sigma) \mp \cos(\alpha\sigma) )( \cosh (2k\sigma) - \cos (2(1-\alpha)\sigma))  },
\end{equation} 
where 
\begin{equation*}
\psi_1 = (k^2+\alpha^2)\sinh (k\sigma) \left[ \sinh (2k\sigma) \mp 2\sinh (k\sigma)\cos(\sigma)\cos((1-\alpha)\sigma)\right]
\end{equation*}
and 
\begin{equation*}
\psi_2 = \frac{k(k^2+\alpha^2)\sinh (k\sigma)}{1-\alpha}[ \sin (2(1-\alpha)\sigma) \mp 2\cosh(k\sigma)\sin((1-\alpha)\sigma)\cos(\sigma) ].
\end{equation*}
From ($\ref{discat0}$) it is clear that for fixed $\alpha\neq 1$
and fixed $\sigma\in(0,2\pi)$, $\bpk$ and $\bmk$ as functions of
$k$ are continuous everywhere except $k=0$.  If we take the limit of ($\ref{discat0}$) as $k$ approaches $0$, we find that
\begin{equation}\label{blowup}
\lim_{k\rightarrow 0} \bk{\pm,k} = \bk{\pm,0}= \frac {\alpha^2\psi_3 } { 2(1\mp \cos(\alpha\sigma) )( 1-\cos 2(1-\alpha)\sigma )},
\end{equation}
where 
\begin{equation}\label{blowup2}
\begin{split}
&\psi_3 = 2\sigma^2 \mp\sigma^2(\cos (\alpha\sigma) + \cos(2-\alpha)\sigma )\\
 &\qquad- \frac{\sigma}{1-\alpha}(\sin(2-2\alpha)\sigma \mp (\sin(2-\alpha)\sigma -\sin(\alpha\sigma))  ).
\end{split}
\end{equation}
From ($\ref{blowup}$) we see that $\beta_{\sigma,\alpha}$ typically blows up when either $\alpha\sigma=n\pi$ or $(1-\alpha)\sigma=n\pi$ for some $n\in\Z$.

The first set of singularities above is a result of the unboundedness of the
Neumann problem for the Laplace operator in weighted spaces on a cone.  To see
this, we observe that in ($\ref{solvefor+and-}$), $\alpha_+$ and $\alpha_-$ become
undefined as $k\to0$ when $\omega^2=e^{-2(k+i\alpha)\sigma}\to1$, which occurs precisely
when $\alpha\sigma=n\pi$ for $n\in\Z$.     

The second set of singularities above, which occur when $(1-\alpha)\sigma=n\pi$, 
result from failure to bound the boundary data $\nddu$ in terms of $\Delta u$.  
For these combinations of $\alpha$ and $\sigma$, the $L^2$ norm of
$r^\alpha\Delta u$ in $\Ks$ is not sufficient to control $\nddu$ appropriately. 
This is fundamentally due to the existence of harmonic fields $u=c
r^{1-\alpha}\sin((1-\alpha)\theta)$ where $c$ is a constant vector.
Corresponding to these fields, there are nontrivial modes $(y_1,y_2)$ for
$k=0$ satisfying
$L_1y_1=0=L_2y_2$ and the no-slip boundary conditions (\ref{BC1}), 
while $(\alpha_+,\alpha_-)$ is non-zero.  
One finds then that the maximum of $I_{p,k}/I_{u,k}\to\infty$ as $k\to0$.

To see just how this can occur in terms of the computations of section 3
for certain combinations of $\sigma$ and $\alpha$ ($\ne1$ or $0$),
we observe that in section 3.5, $L_1y_1=0$ iff $y_1(\theta)$ 
takes the form given in (\ref{26}) with $a_3=a_4=0$.
One can then satisfy the no-slip boundary conditions through (\ref{29}) for
some nonzero $a_1$, $a_2$ if and only if $\omega^2 e^{2i\sigma}=1$, meaning
$k=0$ and $(1-\alpha)\sigma=n\pi$ for some $n\in\Z$.
We may simply take $y_2=0$,
and it follows by (\ref{iuk}) that $I_{u,0}=0$.

But then, $a_1=-a_2\bar\omega\ne0$, and we compute that 
$\partial_{\theta}y_1= 2i(\alpha-1)a_1\ne0$ at $\theta=\sigma$, 
yielding nonzero boundary values for $\partial_\theta\hat q$ in (\ref{partialqhat})
and causing $I_{p,0}$ to be positive in (\ref{Ipk1}).  
If we vary $k$ while holding $(y_1,y_2)$ fixed and use (\ref{yboundary}) to determine
$(\alpha_+,\alpha_-)$ and thence $\hat q$, we see that $I_{u,k}\to0$ as
$k\to0$ while $I_{p,k}\to I_{p,0}>0$.
This results in $\beta_k\to\infty$ as $k\to0$, hence 
$\beta_{\sigma,\alpha}=\infty$ when $(1-\alpha)\sigma=n\pi$.

\section{Proof of Theorem \ref{boundedcase}}\label{bounded}
  
We now use Theorem \ref{main} to prove Theorem \ref{boundedcase}
through a localization argument.  
Let $\Omega$ denote a bounded domain with a straight corner.
Replacing $\Omega$ by a suitable rotated translate if necessary, 
we may assume there is a neighborhood $U$ of $0$ such that 
$U\cap \Omega=U\cap\mathcal{K}_{\sigma}$, where $\sigma\ne\pi$.

Fix any $\beta<1$ and $C\in\R$.  We observe
from the formula for ${\hat{\beta}}_{\pm,k}$ with $\alpha=0$ given in
($\ref{alpha0}$) that $\beta_{\sigma,0}\geq 1$ when $\sigma\neq \pi$.
Therefore, 
there exists a solution ($u, p$) to ($\ref{BVP1}$) with $u$ in
$C^{\infty}_c(\overline{\mathcal{K}}_{\sigma}\backslash \{0\},\R^2)$ which satisfies
$\int_{{\mathcal{K}}_{\sigma}} |\nabla p|^2>\beta \int_{{\mathcal{K}}_{\sigma}}
|\Delta u |^2$.   Replacing $(u,p)$ by suitable dilates if necessary,
we may assume that the support of $u$ is contained in $U$.

We construct a sequence of
solutions ($u_j,p_j$) to ($\ref{BVP1}$) on $\Omega$ by setting 
\begin{equation*}
u_j(x)=j^{-1}u(jx)|_{\Omega} \text{ and }\nabla p_j=(I-P)(\Delta -\nabla\nabla\cdot)u_j\text{ in } \Omega.
\end{equation*}   
We see that $\Delta p_j = 0$ in $\Omega$ and $n\cdot\nabla p_j = n\cdot(\Delta - \nabla\nabla \cdot) u_j$ on $\partial\Omega$.  Moreover, since $u_j$ is supported in $\Omega\cap {\mathcal{K}}_{\sigma}$ for all $j$, we have $\|\Delta u_j\|_{L^2(\Omega)}=\|\Delta u\|_{L^2(j\Omega)}\leq\|\Delta u\|_{L^2({\mathcal{K}}_{\sigma})}$ and $\|\nabla u_j\|_{L^2(\Omega)}= j^{-1}\|\nabla u\|_{L^2(j\Omega)}\leq j^{-1}\|\nabla u\|_{L^2({\mathcal{K}}_{\sigma})}$ for every $j$.  This construction allows us to write the following series of inequalities for sufficiently small $\epsilon>0$ and for sufficiently large $j$: 
\begin{equation*}
\begin{split}
&\int_{\Omega} (\beta |\Delta u_j|^2 + C|\nabla u_j|^2) \leq \int_{{\mathcal{K}}_{\sigma}} (\beta |\Delta u|^2 + Cj^{-2}|\nabla u|^2) \\
&\qquad\qquad \leq \int_{{\mathcal{K}}_{\sigma}} (\beta|\Delta u|^2 )  + \epsilon <  \int_{{\mathcal{K}}_{\sigma}} |\nabla p|^2.
\end{split}
\end{equation*}   
We claim that 
\begin{equation}\label{lowersemicty}
\int_{{\mathcal{K}}_{\sigma}} |\nabla p|^2 \leq \liminf_{j\rightarrow \infty}  \int_{\Omega}|\nabla p_j|^2.
\end{equation}  
To see that ($\ref{lowersemicty}$) holds, we first use the definition of $u_j$, the equality $\nabla p_j = (I-P)(\Delta-\nabla\nabla\cdot)u_j$, and orthogonality of the Leray projection to observe that
\begin{equation}\label{weakconverge}
\int_{\Omega} |\nabla p_j|^2 \leq \int_{j\Omega} |(\Delta-\nabla\nabla\cdot) u|^2 \leq \int_{{\mathcal{K}}_{\sigma}}   |(\Delta-\nabla\nabla\cdot) u|^2
\end{equation}           
for all $j$.  If we define
\begin{equation}
p_j^*(x)= p_j\left(\frac{x}{j}\right) -\frac{1}{m(B)} \int_{B} p_j\left(\frac{x}{j}\right) 
\end{equation}
for $x\in j\Omega$, where $B$ corresponds to the domain $B$ given in (\ref{projunbounded}), then we can apply a generalized Poincare Inequality (see \cite[Ch. 2]{S}) to conclude that  for each $n\in\N$, 
\begin{equation}\label{poincare}
\int_{n\Omega} |p_j^*|^2 \leq C_n \int_{n\Omega} |\nabla p_j^*|^2 \leq C_n \int_{\Omega} |\nabla p_j|^2
\end{equation}
for sufficiently large $j$.  By a standard diagonalization argument, we can construct a subsequence of $\{  \nabla p_j^*\}$, which we henceforth denote as $\{\nabla p_j^* \}$, converging weakly in $L^2(n\Omega)$ for every $n\in\N$.  This implies by ($\ref{poincare}$) and by another diagonalization argument that, up to subsequences, $\{ p^*_j\}$ converges weakly to some $p^*$ in $L^2(n\Omega)$ for all $n\in\N$.  By uniqueness of weak limits, we can conclude that $\{\nabla p_j^*\}$ converges weakly to $\nabla p^*$ in $L^2(n\Omega)$ for every $n$.  Moreover, by properties of weakly convergent sequences we can write $\int_{n\Omega} |\nabla p^{*}|^2 \leq \liminf_{j\rightarrow \infty}  \int_{n\Omega} |\nabla p_j^*|^2$ for each $n$.  We can then conclude that for sufficiently large $j$, 
\begin{equation*}
\begin{split}
\int_{{\mathcal{K}}_{\sigma}} |\nabla p^*|^2 &\leq \lim_{n\rightarrow \infty} \int_{n\Omega} |\nabla p^*|^2 \leq   \liminf_{j\rightarrow \infty} \int_{j\Omega} |\nabla p_j^*|^2=  \liminf_{j\rightarrow \infty} \int_{\Omega} |\nabla p_j|^2. 
\end{split}
\end{equation*}  
It remains to show that $\int_{{\mathcal{K}}_{\sigma}} |\nabla p^{*}|^2=\int_{{\mathcal{K}}_{\sigma}} |\nabla p|^2$.  This will imply ($\ref{lowersemicty}$). 

To show that $\int_{{\mathcal{K}}_{\sigma}} |\nabla p^{*}|^2=\int_{{\mathcal{K}}_{\sigma}} |\nabla p|^2$, we first show $\Delta p^*=0$ in ${\mathcal{K}}_{\sigma}$.  We fix a compact subset $K$ of ${\mathcal{K}}_{\sigma}$, and we apply the mean value property and weak convergence of $\{ p^*_j \}$ to conclude that for any $y\in K$, $|p_j^{*}(y)|\leq C\|p_j^{*}\|_{L^2(n\Omega)}\leq C$, giving equiboundedness of $\{p^{*}_j\}$ on $K$.  Moreover, by the mean value property and weak convergence of $\{ \nabla p^*_j \}$, $\{ \nabla p^*_j \}$ is equibounded on $K$, implying that $\{ p_j^{*} \}$ is also equicontinuous.  Therefore, up to subsequences, $\{ p_j^{*} \}$ converges uniformly on $K$ to $p^{*}$.  We again apply the mean value property and uniform convergence of $\{ p^{*}_j\}$ on $K$ to conclude that $p^{*}$ is harmonic.  

Since $\Delta p^*=0$ on $n\Omega$, we infer that the sequence $\{ \nabla p_j^* \}$ converges weakly to $\nabla p^*$ in $H$(div,$n\Omega)$, the space of vector fields in $L^2(n\Omega)$ with divergence in $L^2(n\Omega)$.  By the boundedness of the trace operator mapping $H$(div, $n\Omega$) into $H^{-\frac{1}{2}}(\partial(n\Omega))$ (see, for example, \cite{GR} Theorem 2.5), we can conclude that $n\cdot\nabla p_j^*$ converges weakly to $n\cdot\nabla p^*$ in $H^{-\frac{1}{2}}(\partial(n\Omega))$.  As $n\cdot\nabla p_j^*=n\cdot\nabla p$ on $\partial{\mathcal{K}}_{\sigma}\cap\partial (n\Omega)$ for every $n$, it follows that $n\cdot\nabla p^*=n\cdot\nabla p$ on $\partial {\mathcal{K}}_{\sigma}$.  

Using the equalities $\Delta p=\Delta p^*=0$ in ${\mathcal{K}}_{\sigma}$ and $n\cdot\nabla (p^*- p)=0$ on $\partial{\mathcal{K}}_{\sigma}$, we can now integrate by parts to conclude that $\int_{{\mathcal{K}}_{\sigma}} |\nabla(p^* - p)|^2=0$.  For $\phi\in C^{\infty}_c(\overline{{\mathcal{K}}}_{\sigma})$, we have that
\begin{equation*}
\int_{{\mathcal{K}}_{\sigma}} \nabla\phi\cdot\nabla (p^*-p) = \int_{\partial{\mathcal{K}}_{\sigma}} \phi n\cdot\nabla (p^*-p) - \int_{{\mathcal{K}}_{\sigma}} \phi\Delta(p^*-p)=0.
\end{equation*}        
Since $p^*-p$ belongs to $Y$ and $C^{\infty}_c(\overline{{\mathcal{K}}}_{\sigma})$ is dense in $Y$, it follows that $\int_{{\mathcal{K}}_{\sigma}} |\nabla (p^*-p)|^2=0$, and ($\ref{lowersemicty}$) holds.  This completes the proof of Theorem \ref{boundedcase}.
          

\bibliographystyle{amsplain}

\end{document}